\documentclass[11pt,a4paper,reqno]{amsart}
\usepackage[belowskip=-15pt,aboveskip=0pt]{caption}
\usepackage{graphicx}
\usepackage{amsmath}
\usepackage{amsfonts}
\usepackage{amssymb}
\usepackage{amscd}
\usepackage{url}
\usepackage[pdftex,bookmarks=true]{hyperref}
\usepackage{enumerate}

\newcommand\vol{{\operatorname{vol}}}

\newcommand\R{{\mathbf{R}}}
\newcommand\C{{\mathbf{C}}}

\renewcommand\P{{\mathbf{P}}}
\newcommand\E{{\mathbf{E}}}

\newcommand\tr{{\operatorname{tr}}}

\newcommand\dist{{\operatorname{dist}}}
\newcommand\Z{{\mathbf{Z}}}

\newcommand\col{{\mathbf{c}}}
\newcommand\row{{\mathbf{r}}}
\newcommand\Var{{\mathbf{Var}}}

\newcommand\Ba{{\mathbf a}}
\newcommand\Bb{{\mathbf b}}

\newcommand\Be{{\mathbf e}}
\newcommand\Bf{{\mathbf f}}

\newcommand\Bk{{\mathbf k}}

\newcommand\Br{{\mathbf r}}

\newcommand\Bu{{\mathbf u}}
\newcommand\Bv{{\mathbf v}}

\newcommand\Bx{{\mathbf x}}
\newcommand\By{{\mathbf y}}

\newcommand\BG{{\mathbf G}}
\newcommand\BH{{\mathbf H}}

\newcommand\cir{{\mathbf{cir}}}

\newcommand\ep{\epsilon}
\newcommand\Tr{{\mathbf{Tr}}}
\theoremstyle{plain}
  \newtheorem{theorem}[subsection]{Theorem}

  \newtheorem{fact}[subsection]{Fact}
  \newtheorem{lemma}[subsection]{Lemma}
  \newtheorem{corollary}[subsection]{Corollary}
  
  \newtheorem{example}[subsection]{Example}

\theoremstyle{remark}
  \newtheorem{remark}[subsection]{\bf Remark}
  
  \newtheorem{claim}[subsection]{\bf Claim}

\theoremstyle{definition}

\include{psfig}

\begin{document}

\title[random matrices of given row sum]{Circular law for random discrete matrices of given row sum}

\author{Hoi H. Nguyen}
\address{Department of Mathematics, University of Pennsylvania, Philadelphia, PA 19104, USA}

\email{hoing@math.upenn.edu}
\author{Van H. Vu}
\address{Department of Mathematics, Yale University, New Haven , CT 06520, USA}
\email{van.vu@yale.edu}

\maketitle

\begin{abstract} Let $M_n$ be a random matrix of size $n\times n$ and let $\lambda_1,\dots,\lambda_n$ be the eigenvalues of $M_n$. The {\it empirical spectral distribution} $\mu_{M_n}$ of $M_n$ is defined as

$$\mu_{M_n}(s,t)=\frac{1}{n}\# \{k\le n, \Re(\lambda_k)\le s; \Im(\lambda_k)\le t\}.$$ 

The circular law theorem in random matrix theory asserts that if the entries of $M_n$ are i.i.d. copies of a random variable with mean zero and variance $\sigma^2$, then the empirical spectral distribution of the normalized matrix $\frac{1}{\sigma\sqrt{n}}M_n$ of $M_n$ converges almost surely to the uniform distribution $\mu_\cir$ over the unit disk as $n$ tends to infinity. 

\vskip .1in

In this paper we show that the empirical spectral distribution of the normalized matrix of $M_n$, a random matrix whose rows are independent random $(-1,1)$ vectors of given row-sum $s$ with some fixed integer $s$ satisfying $|s|\le (1-o(1))n$,  also obeys the circular law. The key ingredient is a new polynomial estimate on the least singular value of $M_n$.  
\end{abstract}

\section{Introduction}\label{section:introduction}

Let $M_n$ be a matrix of size $n\times n$ and let $\lambda_1,\dots,\lambda_n$ be the eigenvalues of $M_n$. Then the empirical spectral distribution (ESD) $\mu_{M_n}$ of $M_n$ is defined as

$$\mu_{M_n}(s,t)=\frac{1}{n}\# \{k\le n, \Re(\lambda_k)\le s; \Im(\lambda_k)\le t\}.$$ 

We also define $\mu_{\cir}$ as the uniform distribution over the unit disk,

$$\mu_{\cir}(s,t)= \frac{1}{\pi} mes(|z|\le 1; \Re(z)\le s, \Im(z)\le t).$$ 

Confirming a long standing conjecture in random matrix theory, a recent result of Tao and Vu (appendix by Krishnapur) proves a universal law for the ESD of random i.i.d. matrices.

\begin{theorem}\cite{TVcir}\label{theorem:TVcir}
Assume that the entries of $M_n$ are i.i.d. copies of a complex random variable of mean zero and finite non-zero variance $\sigma^2$, then the ESD of the matrix $\frac{1}{\sigma \sqrt{n}}M_n$ converges to $\mu_\cir$ almost surely as $n$ tends to $\infty$.
\end{theorem}

The proof of this result is built upon previous important developments of Girko \cite{G1,G2}, Bai \cite{B}, G\"otze-Tikhomirov \cite{GT}, Pan-Zhou \cite{PZ}, Tao-Vu \cite{TVcir0} and many others. 

In view of universality phenomenon, it is of importance to study the law for random matrices of non-independent entries. Probably one of the first results in this direction is due to Bordenave, Caputo and Chafai \cite{BCC} who prove the law for random Markov matrices.

\begin{theorem}\cite[Theorem 1.3]{BCC} \label{theorem:BCC}
Let $X$ be a random matrix of size $n \times n$ whose entries are i.i.d. copies of a non-negative continuous random variable with finite variance $\sigma^2$ and bounded density function. Then with probability one the ESD of the normalized matrix $\sqrt{n}\bar{X}$, where $\bar{X}=(\bar{x}_{ij})_{1\le i,j \le n}$ and ${\bar{x}}_{ij}:=x_{ij}/(x_{i1}+\dots+x_{in})$, converges weakly to the circular measure $\mu_\cir$.
\end{theorem}

In particular, when $x_{11}$ follows the exponential law of mean one, Theorem \ref{theorem:BCC} establishes the circular law for the Dirichlet Markov ensemble (see also \cite{C}). We remark that the assumptions of continuity and boundedness are crucial in the proof of Theorem \ref{theorem:BCC}.

Related results with ''linear'' assumption of independence include a result of Tao, who among other things proves the circular law for random zero-sum matrices.

\begin{theorem}\cite[Theorem 1.13]{Tao}   
Let $X$ be a random matrix of size $n \times n$ whose entries are i.i.d. copies of a random variable of mean zero and variance one. Then the ESD of the normalized matrix $\frac{1}{\sqrt{n}}\bar{X}$, where $\bar{X}=(\bar{x}_{ij})_{1\le i,j \le n}$ and $\bar{x}_{ij}:=x_{ij}-\frac{1}{n}(x_{i1}+\dots+x_{in})$, converges almost surely to the circular measure $\mu_\cir$.
\end{theorem}

The main goal of this note is to showing that the circular law also holds for random discrete matrices of similar weak constraints.

\begin{theorem}[Main result]\label{theorem:cir} Let $0<\ep\le 1$ be a positive constant. Let $M_n$ be a random $(-1,1)$ matrix of size $n\times n$ whose rows are independent vectors of given row-sum $s$ with some $s$ satisfying $|s|\le (1-\ep)n$. Then the ESD of the normalized matrix $\frac{1}{\sigma\sqrt{n}}M_n$, where $\sigma^2=1-(\frac{s}{n})^2$, converges almost surely to the distribution $\mu_{\cir}$ as $n$ tends to $\infty$.
\end{theorem}

To some extent, our matrix is a discrete version of the random Markov matrices considered in Theorem \ref{theorem:BCC} where the entries are restricted to $\pm 1/s$. However, it is probably more suitable to compare our model with that of random Bernoulli matrices. By Theorem \ref{theorem:TVcir}, the ESD of the normalized random Bernoulli matrices obeys the circular law, and hence our Theorem \ref{theorem:cir} serves as a local version of the law.

We remark that in a very recent result \cite{Ng-doubly}, the first author is able to prove a similar law for random doubly stochastic matrices, thus confirming the universality principle for another type of matrix of independent entries. Although the results are similar in spirit, the difficulties in each note are very different. The main obstacle of this note is to study the singularity of $M_n$ and its perturbed variants. Inverse techniques developed in the literature to deal with this problem do not seem to suffice. This leads us to a new development to be discussed in Section \ref{section:singular}. Note that our approach may also cover the regime $n-s=o(n)$ but we do not attempt to do so here. In what follows we present some reduction steps to simplify our problem.

\begin{figure}[!ht]
\centering
\includegraphics[scale=0.5]{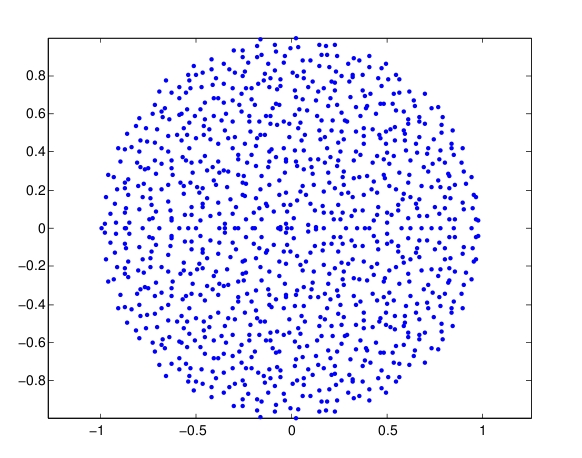}
\vskip .2in
\caption{The ESD of a random matrix of size 1000 by 1000 whose rows are $(-1,1)$ vectors of zero-sum, picture by Phillip Woods.}
\end{figure}

Observe that, by letting $X_{n-1}$ be the submatrix generated by the first $n-1$ rows and columns of $M_n$, the spectra of $M_n$ is the union of $s$ and the spectra of the pertubed matrix $X_{n-1}-F_{n-1}$ where all of the rows of $F_{n-1}$ are identical copies of $(m_{n1},\dots,m_{n(n-1)})$, here by $m_{ij}$ we mean the $ij$-th entry of $M_n$. 

Indeed, consider the matrix $M:=M_n - \lambda I_n$. We have 

$$\det(M)=\det(M'),$$ 

where $M'$ is obtained from $M$ by adding its first $n-1$ columns to its last one. 

On the other hand, we also have

$$\det(M')=(s-\lambda)\det(M''),$$ 

where 

$$ M'':=
\begin{pmatrix}
  m_{11}-\lambda & \cdots & m_{1(n-1)} & 1 \\
  \vdots  & \ddots  & \vdots & \vdots  \\
  m_{(n-1)1} & \cdots & m_{(n-1)(n-1)}-\lambda & 1 \\
  m_{n1} & \cdots & m_{n(n-1)} & 1 \\
\end{pmatrix}.
$$

It is  clear that $\det(M'') = \det(M''')$,  where $M''':=(X_{n-1}-F_{n-1})-\lambda I_{n-1}$. Thus the spectra of $M_n$ is indeed the union of $s$ and the spectra of the pertubed matrix $X_{n-1}-F_{n-1}$.

The observation above suggests a way to prove Theorem \ref{theorem:cir} by looking at the ESD of $X_{n-1}-F_{n-1}$. This alternative helps us avoid the outlier eigenvalue $s$ of $M_n$ which may cause certain technical difficulty for any direct study on $M_n$. 

Notice that the rows of $X_{n-1}$ above are independent vectors chosen uniformly from the set of all $(-1,1)$ vectors of row-sum either $s-1$ or $s+1$. So for Theorem \ref{theorem:cir} it suffices to show the following.

\begin{theorem}[Circular law for pertubed matrices]\label{theorem:cir:X} Let $X_n$ be a random $(-1,1)$ matrix whose rows are independent random vectors of row-sum either $s-1$ or $s+1$ with given $s$ satisfying $|s|\le (1-\ep)n$. Let $F_n$ be a deterministic matrix whose rows are identical copies of a given $(-1,1)$ vector $\Bf$. Then the ESD of $\frac{1}{\sigma\sqrt{n}}(X_n+F_n)$, where $\sigma^2=1-(\frac{s}{n})^2$, converges almost surely to the distribution of $\mu_{\cir}$ as $n$ tends to $\infty$.
\end{theorem}

For short, by $\mathcal{S}$ we denote the set of all $(-1,1)$ vectors $\Bx=(x_1,\dots,x_n)$ of row-sum either $s-1$ or $s+1$. To establish Theorem \ref{theorem:cir:X} we will relate $X_n$ to a random matrix $X_n'$ whose entries are i.i.d. copies of a random Bernoulli variable $x$ of the following form

\begin{equation}\label{eqn:Bernoulli}
\begin{cases}
\P(x=-1) = \frac{1}{2}-\frac{s}{2n},\\
\P(x=1)=\frac{1}{2}+\frac{s}{2n}.
\end{cases}
\end{equation}

It is known that the ESD of $\frac{1}{\sigma\sqrt{n}}(X_n'+F_n)$ converges uniformly to $\mu_\cir$ (see for instance \cite[Corollary 1.15]{TVcir}). As we desire to pass this result to $X_n+F_n$, we will make use of a so called replacement principle below. 
 
\begin{theorem}\cite[Theorem 2.1]{TVcir}\label{theorem:replacement}
Suppose for each $n$ that $A_n=(a_{ij}),B_n=(b_{ij})$ are random matrices of size $n\times n$. Assume that 

\begin{itemize}
\item the sum

$$\frac{1}{n^2}\sum_{ij}(|a_{ij}|^2+ |b_{ij}|^2)$$

\noindent is bounded almost surely;

\vskip .1in

\item for almost all complex numbers $z$

$$\frac{1}{n}\log |\det(\frac{1}{\sqrt{n}}A_n-zI_n)| - \frac{1}{n}\log |\det(\frac{1}{\sqrt{n}}B_n-zI_n)|$$

\noindent converges almost surely to zero.
\end{itemize}

Then $\mu_{\frac{1}{\sqrt{n}}A_n} - \mu_{\frac{1}{\sqrt{n}}B_n} $ converges almost surely to zero.
\end{theorem}

In application, $X_n+F_n$ plays the role of $A_n$ and $X_n'+F_n$ plays that of $B_n$. It is clear that the first condition of Theorem \ref{theorem:replacement} is satisfied. Thus for Theorem \ref{theorem:cir:X} it suffices to justify the second condition.

\begin{theorem}\label{theorem:comparison} For every fixed complex $z$ we have

$$\frac{1}{n}\log |\det((X_n+F_n)-z\sqrt{n}I_n)| - \frac{1}{n}\log |\det((X_n'+F_n)-z\sqrt{n}I_n)|$$

converges to zero almost surely.
\end{theorem}

We will outline a proof for Theorem \ref{theorem:comparison} in the next section.

{\bf Notation.} Here and later, asymptotic notations such as $O,\Omega,\Theta$, and so for, are used under the assumption that $n\rightarrow \infty$. A notation such as $O_C(.)$ emphasizes that the hidden constant in $O$ depends on $C$. 

For $1\le s \le n$, we denote by $\Be_s$ the unit vector $(0,\dots,0,1,0,\dots,0)$, where all but the $s$-th component are zero. For a real or complex vector $\Bv=(v_1,\dots,v_n)$, we use the shorthand $\|\Bv\|$ for its $L_2$-norm $(\sum_i|v_i|^2)^{1/2}$.  

For a matrix $M$, we use the notation $\row_i(M)$ and $\col_j(M)$ to denote its $i$-th row and $j$-th column respectively. For an event $A$, we use the subscript $\P_{\Bx}(A)$ to emphasize that the probability under consideration is taking according to the random vector $\Bx$.

\section{Proof of Theorem \ref{theorem:comparison}: outline}\label{section:comparison}

Let $\Bf_1,\dots,\Bf_n$ denote the (deterministic) rows of $F_n+\sqrt{n}zI_n$, and let $\Bx_1,\dots,\Bx_n$ as well as $\Bx_1',\dots,\Bx_n'$ be the rows of $X_n$ and $X_n'$ respectively.

For each $i\ge 2$, let $V_{i-1}$  be the space spanned by $\Bx_1+\Bf_1,\dots,\Bx_{i-1}+\Bf_{i-1}$ and let $\dist(\Bx_i+\Bf_i,V_{i-1})$ be the distance from $\Bx_i+\Bf_i$ to $V_{i-1}$. Define similarly for $V_{i-1}'$ and $\dist(\Bx_i'+\Bf_i,V_{i-1}')$.  By the ''base times height'' formula we have

\begin{align*}
&\log \Big|\det((X_n+F_n)-z\sqrt{n}I_n)\Big| = \sum_{i} \log \dist((\Bx_i+\Bf_i),V_{i-1}).\\
&=\sum_{i\le m} \log \dist((\Bx_i+\Bf_i),V_{i-1}) + \sum_{m< i} \log \dist((\Bx_i+\Bf_i),V_{i-1})\\
&:= \log S_1+ \log S_2;
\end{align*}

and similarly,

\begin{align*}
&\log \Big|\det((X_n'+F_n)-z\sqrt{n}I_n)\Big| = \sum_{i} \log \dist((\Bx_i'+\Bf_i),V_{i-1}).\\
&=\sum_{i\le m} \log \dist((\Bx_i'+\Bf_i),V_{i-1}') + \sum_{m< i} \log \dist((\Bx_i'+\Bf_i),V_{i-1}')\\
&:= \log S_1'+ \log S_2'.
\end{align*}

where we set the threshold $m$ to be $m:=n-\log^8 n$.

In order to compare $\log \Big|\det((X_n+F_n)-z\sqrt{n}I_n)\Big|$ with $\log \Big|\det((X_n'+F_n)-z\sqrt{n}I_n)\Big|$ we will show the following.

\begin{theorem}\label{theorem:large} With probability $1-\exp(-\log^{2-o(1)}n)$ we have

$$\frac{1}{n} |\log S_1 - \log S_1'|=O(\log^{-2}n).$$
\end{theorem}

\begin{theorem}\label{theorem:small} With probability  $1-O(n^{-100})$ we have
$$\frac{1}{n}(|\log S_2| + |\log S_2'|)=O(\log^9 n/n).$$ 
\end{theorem}

It is clear that Theorem \ref{theorem:comparison} follows from Theorem \ref{theorem:large} and Theorem \ref{theorem:small}. In what follows we outline the approach to prove these results.

\subsection{Sketch of the proof of Theorem \ref{theorem:large}} One of the main ingredients is the following row replacement principle.

\begin{lemma}\label{lemma:replacement} Let $i$ be an integer between $1$ and $m$. Let $\Bx_1,\dots,\Bx_i,\Bx_{i}',\Bx_{i+1}',\dots,\Bx_m'$ be $m+1$ independent vectors where the $\Bx_j$'s are random vectors of type $\mathcal{S}$ and $\Bx_k'$'s are random vectors whose components are i.i.d copies of $x$ from \eqref{eqn:Bernoulli}. Assume that $\vol_i$ is the $m$-dimensional volume of the parallelepiped generated by $\Bx_1+\Bf_1,\dots,\Bx_i+\Bf_i, \Bx_{i+1}'+\Bf_{i+1},\dots,\Bx_m'+\Bf_m$ and $\vol_{i-1}$ is that of the parallelepiped generated by $\Bx_1+\Bf_1,\dots,\Bx_{i-1}+\Bf_{i-1},\Bx_i'+\Bf_i,\dots,\Bx_m'+\Bf_m$. Then we have

$$\P_{\Bx_1,\dots,\Bx_i,\Bx_{i}',\Bx_{i+1}',\dots,\Bx_m'}\Big(|\log \vol_i -\log \vol_{i-1}|=O(\log^{-2}n)\Big) = 1-\exp(-\log^{2-o(1)}n).$$
\end{lemma}

Lemma \ref{theorem:large} then follows by a repeatedly use of Lemma \ref{lemma:replacement} and the triangle inequality using the fact that $S_1$ and $S_1'$ are volumes of the parallelepipeds generated by  $\Bx_1+\Bf_1,\dots,\Bx_m+\Bf_m$ and by $\Bx_1'+\Bf_1,\dots,\Bx_m'+\Bf_m$ respectively.

We now justify Lemma \ref{lemma:replacement}. We express $\vol_i$ as $\vol_i=d \times \vol$, where $d$ is the distance from $\Bx_{i}+\Bf_{i}$ to the space $V$ spanned by $\Bx_1+\Bf_1,\dots,\Bx_{i-1}+\Bf_{i-1},\Bx_{i+1}'+\Bf_{i+1}',\dots,\Bx_m'+\Bf_m$ and $\vol$ is the volume of the parallelepiped generated by these vectors. Similarly we can express $\vol_{i-1}$ as $\vol_{i-1}=d'\times \vol$, where $d'$ is the distance from $\Bx_i'+\Bf_i$ to $V$. 
 
Thus we have 

$$|\log (\vol_i) -\log (\vol_{i-1})|=|\log d -\log d'|.$$ 

We will next see that $d$ and $d'$ are almost identical with very high probability. 

Let $\Bf$ be a fixed vector (whose coordinates may depend on $n$). In what follows we denote the translation $\Bf+(s/n,\dots,s/n)$ of $\Bf$ by $\Bf'$.

\begin{lemma}\label{lemma:Talagrand} Assume that $V\subset \C^n$ is a subspace of dimension $\dim(V)=k \le n-10$. Let $\Bx'=(x_1',\dots,x_n')$ be a random vector where $x_i'$ are i.i.d. copies of $x$ from \eqref{eqn:Bernoulli} and let $d'$ be the distance from $\Bx'+\Bf$ to $V$. Then for any $t>0$ we have

$$\P_{\Bx'}(|d'-\sqrt{n-k + d_{\Bf'}^2}|\ge t+3) \le \exp(-\frac{t^2}{4}),$$

where $d_{\Bf'}$ is the distance from $\Bf'$ to $V$.

\end{lemma}

Lemma \ref{lemma:Talagrand} can be proved by using a well-known result of Talagrand; we defer its proof to Section \ref{section:Talagrand}. 

As $\E(\sum_ix_i')=s$ and $\Var(\sum_ix_i')=\Theta(n)$, the probability that a random vector $\Bx'$ belongs to the set of $(-1,1)$ vectors of row-sum $s+1$ (or $s-1$) is $\Theta(1/\sqrt{n})$. Furthermore, condition on $\Bx' \in \mathcal{S}$, $\Bx'$ is uniformly distributed over these sets. We thus infer from Lemma \ref{lemma:Talagrand} the following.  

\begin{corollary}\label{cor:Talagrand}
Let $\Bx$ be a vector uniformly sampled from $\mathcal{S}$ and let $d$ be the distance from $\Bx+\Bf$ to $V$. Then for any $t>0$ we have

$$\P_{\Bx}(|d-\sqrt{n-k+d_{\Bf'}^2}|\ge t+3) =O(\sqrt{n} \exp(-\frac{t^2}{4})).$$
\end{corollary}

One immediate consequence of Lemma \ref{lemma:Talagrand} and Corollary \ref{cor:Talagrand} is that if $k\le n-\log^4 n$, then by setting $t=\log n$, $d$ is nonzero with probability at least $1-O(\exp(-\log^{2-o(1)}n))$. By applying this fact $m$ times, we conclude that all the $\vol_i$ are non-zero with probability at least $1-O(\exp(-\log^{2-o(1)}n)$. So it is safe to assume that $V$ has dimension exactly $m-1$ for any $V$ spanned by $\Bx_1+\Bf_1,\dots,\Bx_{i-1}+\Bf_{i-1},\Bx_{i+1}'+\Bf_{i+1}',\dots,\Bx_m'+\Bf_m$. Next, by applying Lemma \ref{lemma:Talagrand} and Corollary \ref{cor:Talagrand} once more,  with probability $1-O(\exp(-\log^{2-o(1)}n))$ with respect to $\Bx_i$ and $\Bx_i'$ we have 

$$|d-\sqrt{n-m+1+d_{\Bf_i'}^2}|\le\log n$$ 

and 

$$|d'-\sqrt{n-m+1+d_{\Bf_i'}^2}|\le\log n.$$

It then follows that 

$$|\log d -\log d'| \le \log (1+\frac{2\log n}{\log^4n -\log n}) = O(\log^{-2} n),$$
 
completing the proof of Lemma \ref{lemma:replacement}.

\subsection{Sketch of the proof of Theorem \ref{theorem:small}} Our key lemma here is to showing that the least singular value of $X_n+F_n+z\sqrt{n}I_n$, for any fixed complex number $z$, is at least $n^{-O(1)}$ with probability $1-O(n^{-100})$. 

\begin{theorem}\label{theorem:singular'} Assume that $F$ is a deterministic complex  matrix of size $n\times n$ such that $|f_{ij}|\le n^{\gamma}$ for some constant $\gamma$. Then for any $B>0$ there exists $A>0$ depending on $B$ and $\gamma$ such that 

$$\P\big(\sigma_n(X_n+F)<n^{-A}\big)\le O(n^{-B}).$$ 
\end{theorem}

This theorem is an analog of the Bernoulli counterpart $X_n'+F$ whose proof can be found in either \cite{TVcomp} or in other papers of the second author with Tao such as \cite{TVbull,TVcir,TVinverse}. Unfortunately, these proofs do not seem to cover Theorem \ref{theorem:singular} in any trivial way. Henceforth a large part of this note will be devoted to prove it, starting from Section \ref{section:singular}.

%Although bounding the least singular values is not a straightforward task, the largest singular values $\sigma_1(X_n+F_n+z\sqrt{n}I_n)$ and  $\sigma_1(X_n'+F_n+z\sqrt{n}I_n)$ can be easily bounded by $n^{O(1)}$ as the entries of the matrices are all bounded by $O(\sqrt{n})$.

We next invoke the following two linear algebra results.

\begin{lemma}[Cauchy's interlacing law]\cite[Lemma A.1]{TVcir}\label{lemma:Cauchy}
Let $A$ be a matrix of size $n\times n$  and $A'$ be the submatrix formed by the first $n-k$ rows of $A$. Let $\sigma_1(A)\ge \dots \ge \sigma_n(A)\ge 0$ be the singular values of $A$, and similarly for $A'$. Then we have 

$$\sigma_i(A)\ge \sigma_i(A') \ge \sigma_{i+k}(A)$$

for every $1\le i \le n-k$. 
\end{lemma}

\begin{lemma}[Negative second moment]\cite[Lemma A.4]{TVcir}\label{lemma:negative}
Let $1\le n' \le n$, and let $A'$ be a full rank matrix of size $n'$ by $n$ with singular values $\sigma_1(A') \ge \dots \ge \sigma_n(A')\ge 0$ and rows $\Br_1,\dots,\Br_{n'}\in \C^n$. For each $1\le i\le n'$, let $W_i$ be the subspace generated by the $n'-1$ rows $\Br_1,\dots,\Br_{i-1},\Br_{i+1},\dots, \Br_{n'}$. Then we have

$$\sum_{i=1}^{n'}\sigma_i^{-2}(A') = \sum_{i=1}^{n'}\dist^{-2}(\Br_i,W_i).$$  
\end{lemma}

We now prove Theorem \ref{theorem:small}. By Theorem \ref{theorem:singular'} we can assume that $\Bx_1+\Bf,\dots,\Bx_n+\Bf$ spans the whole space $\R^n$ with probability at least $1-O(n^{-100})$, and so in particular all the $V_i$ have full rank. Applying Lemma \ref{lemma:negative} for the matrix $A'$ generated by the first $k$ rows $\Bx_1+\Bf,\dots,\Bx_k+\Bf$ with any $k>m=n-
\log^8 n$, we obtain the following with probability at least $1-O(n^{-100})$ 

$$\dist^{-2}(\Bx_k+\Bf,V_{k-1}) < \sum_{i=1}^{k}\sigma_i^{-2}(A')= O(n^{O(1)}),$$

where in the RHS estimate we applied Lemma \ref{lemma:Cauchy} and then Theorem \ref{theorem:singular'}.

Thus for any $k>m$

\begin{equation}\label{eqn:dist}
O(n^{-O(1)})=\dist(\Bx_k+\Bf,V_{k-1}) \le \|\Bx_k+\Bf\|=O(\sqrt{n}).
\end{equation}

Similarly, by applying the known variant of Theorem \ref{theorem:singular'} for $(X_n'+F_n)-z\sqrt{n}I_n$ and by Lemmas \ref{lemma:Cauchy} and \ref{lemma:negative} we also have

\begin{equation}\label{eqn:dist'}
O(n^{-O(1)})=\dist(\Bx_k'+\Bf,V_{k-1}')=O(\sqrt{n}).
\end{equation}

Owing to the estimates \eqref{eqn:dist} and \eqref{eqn:dist'}, we infer that 

$$\P\Big(\frac{1}{n}(|\log S_2| + |\log S_2'|)=O(\log^9 n/n)\Big) = 1-O(n^{-100}),$$ 

proving Lemma \ref{theorem:small}.

\section{The least singular value bound}\label{section:singular}
For the reader's convenience, we restate Theorem \ref{theorem:singular'} below.

\begin{theorem}\label{theorem:singular} Assume that $F$ is a deterministic complex matrix such that $|f_{ij}|\le n^{\gamma}$ for some constant $\gamma$. Then for any $B>0$ there exists $A>0$ depending on $B$ and $\gamma$ such that 

$$\P\big(\sigma_n(X_n+F)<n^{-A}\big)\le O(n^{-B}).$$ 
\end{theorem}

This section is devoted to provide an overview of our approach to prove Theorem \ref{theorem:singular}. More details of the proofs will be discussed in subsequent sections.

We use the shorthand $X$ for the matrix $X_n+F$. To prove Theorem \ref{theorem:singular}, we assume that there exist vectors $\Ba$ and $\Bb$ in $\C^n$ such that $\|\Ba\|=1, \|\Bb\| <  n^{-A}$ and

$$X\Ba=\Bb.$$ 

We next consider two cases.

{\bf Case 1}. $X$ is non-singular. Let $C(X)=(c_{ij}(X))$, $1\le i,j\le n$, be the matrix of the cofactors of $X$. We then have

$$C(X)\Bb = \det(X) \cdot \Ba.$$ 

Thus 

$$\|C(X)\Bb\| = |\det(X)|.$$ 

By paying a factor of $n$ in probability, without loss of generality we can assume that 

$$|c_{11}(X)b_1+\dots c_{1n}(X)b_n|\ge |\det(X)|/n^{1/2}.$$

Note that $\|\Bb\|\le n^{-A}$, thus by Cauchy-Schwarz inequality

\begin{equation}\label{eqn:intro:1}
\sum_{i=1}^n |c_{1i}(X)|^2 \ge n^{2A-1} \det(X)^2.
\end{equation}

We next express $\det(X)$ as a linear form of its first row $\row_1(X)=(x_1+f_{11},\dots,x_n+f_{1n})$

$$\det(Q) = x_1c_{11}(X)+\dots+x_nc_{1n}(X)+r_{11}c_{11}(X)+\dots+r_{1n}c_{1n}(X).$$

Thus, with $c:=\sqrt{\sum_{j}c_{1i}(X)^2}$  (which is $\neq 0$ as  $(c_{11},\dots,c_{1n}) \neq \mathbf{0}$), \eqref{eqn:intro:1} can be rewritten as

$$\Big|x_1 \frac{c_{11}(X)}{c}+\dots+x_n \frac{c_{1n}(X)}{c} + \frac{1}{c}(f_{11}c_{11}(X)+\dots+f_{1n}c_{1n}(Q))\Big| \le n^{-A+1/2}.$$

Roughly speaking, our approach to prove Theorem \ref{theorem:singular} consists of two main steps.
  
\begin{itemize}

\item {\it Step 1}. Condition on $X'$, the matrix of the last $n-1$ rows of $X$, if 

$$\sup_{v}\P_{x_{1},\dots,x_{n}}\big(|\sum_{i=1}^n x_i  \frac{c_{1i}(X_n)}{c} -v|\le  n^{-A}\big)\ge n^{-B},$$ 

then there is a strong structure among the cofactors $c_{1i}$.

\vskip .2in

\item {\it Step 2}. The probability, with respect to $X'$, that there is a strong additive structure among the $c_{1i}$ is negligible.
\end{itemize}

We pause to discuss the structure mentioned in the inverse step. A set $Q\subset \C$ is a \emph{GAP of rank $r$} if it can be expressed as in the form
$$Q= \{g_0+ k_1g_1 + \dots +k_r g_r| k_i\in \Z, K_i \le k_i \le K_i' \hbox{ for all } 1 \leq i \leq r\}$$ 
for some $(g_0,\ldots,g_r) \in \C^{r+1}$ and $(K_1,\ldots,K_r), (K'_1,\ldots,K'_r) \in \Z^r$.

It is convenient to think of $Q$ as the image of an integer box $B:= \{(k_1, \dots, k_r) \in \Z^r| K_i \le k_i \le K_i' \} $ under the linear map $\Phi: (k_1,\dots, k_r) \mapsto g_0+ k_1g_1 + \dots + k_r g_r$.

The numbers $g_i$ are the \emph{generators } of $Q$, the numbers $K_i'$ and $K_i$ are the \emph{dimensions} of $Q$, and $\vol(Q) := |B|$ is the \emph{size} of $B$. We say that $Q$ is \emph{proper} if this map is one to one, or equivalently if $|Q| = \vol(Q)$.  For non-proper GAPs, we of course have $|Q| < \vol(Q)$. If $-K_i=K_i'$ for all $i\ge 1$ and $g_0=0$, we say that $Q$ is {\it symmetric}.

We are now ready to state our steps in details.

\begin{theorem}[Step 1]\label{theorem:step1}
Let $0<\alpha<1/2$ be a given constant. Assume that

$$\rho_{n^{-A}}^\ast(\{v_1,\dots,v_n\}):= \sup_v \P_{x_{1},\dots,x_{n}}\big(|\sum_{i=1}^n x_i v_i-v|\le n^{-A}\big)\ge n^{-B}$$ 

for some sufficiently large $A$, where $v_{i}=c_{1i}(X)/c$. Then, there exists a vector $\Bu=(u_1,\dots,u_{n})$ and a real number $\beta$ of the form $n^{-A+k(5B+5+\gamma)}$ where $0\le k \le A/(10B+10+2\gamma), k\in \Z$ such that the following holds.

\begin{itemize}
\item $\|\Bu\|\asymp 1$ and $|\langle \Bu,\row_i(X)\rangle| \le \beta n^{5B+4+\gamma} $ for $n-1$ rows $\row_i$ of $X$.
\vskip .2in
\item There exists a generalized arithmetic progression $Q^\ast$ of rank $O_{\alpha,B}(1)$ and size  $|Q^\ast|=\max\Big(1,O_{\alpha,B}\big((\rho_{\beta n^{5B+4+\gamma}}^\ast(\{u_1,\dots,u_n\}))^{-1}/n^{\alpha/2}\big)\Big)$ which contains at least $n-n^{1/2+\alpha}$ complex numbers $u_i$.
\vskip .2in
\item All the components of $u_i$ and of the generators of $Q^\ast$ are rational numbers of the form $p/q$, where $|p|,|q| \le n^{A+1}$.
\end{itemize}

\end{theorem}

\vskip .2in

Roughly speaking, the quantity $(\rho_{\beta n^{5B+4+\gamma}}^\ast(\{u_1,\dots,u_n\}))^{-1}$ appearing in the bound of $|Q^\ast|$ guarantees that the containment is economical. 

In the second step of the approach, we  show that the probability for $Q'$ having the above properties is negligible.

\begin{theorem}[Step 2]\label{theorem:step2}
With respect to $X'$, the probability that there exists a vector $\Bu$ and a number $\beta$ as in Theorem \ref{theorem:step1} is $\exp(-\Omega(n))$.
\end{theorem}

We remark here that the choice of $\alpha$ being near $1/2$ would optimize the probability bound in Theorem \ref{theorem:step2}. However, we prefer to keep $\alpha$ abstract to demonstrate the flexibility of our approach.

We now study the remaining case.

{\bf Case 2.} $X$ is singular. We show that the probability of this event is bounded by $O(n^{-B})$ for any $B>0$, where the implied constant depends on $B$. The approach is identical (if not easier) to that of {\bf Case 1}. 

First of all, by paying a factor of $n$ in probability and without loss of generality, it suffices to consider the event that $\Bx_1+\Bf_1$ belongs to the subspace generated by $\Bx_2+\Bf_2,\dots,\Bx_n+\Bf_n$. We show

\begin{theorem}\label{theorem:discrete}
Assume that $X_n$ is a random matrix whose rows $\Bx_1,\dots,\Bx_n$ are independent random vectors sampled uniformly from ${\mathcal S}$. Then for any $B>0$
$$\P(\Bx_1+\Bf_1 \mbox{ belongs to the subspace } H  \mbox{ generated by } \Bx_2+\Bf_2, \dots,\Bx_n+\Bf_n)=O(n^{-B}),$$

where the implied constant depends on $B$.
\end{theorem}

Condition on $\Bx_2,\dots,\Bx_n$, let $\Bv=(v_1,\dots,v_{n})$ be a unit vector which is orthogonal to $H$. Then the probability that $\Bx_{1}+\Bf_1=(x_1+f_{11},\dots,x_n+f_{1n})$ belongs to $\BH$ is bounded by $\P_{x_1,\dots,x_n}(x_1v_1+\dots+x_nv_n + (f_{11}v_1+\dots+f_{1n}v_n)=0)$, and so crudely  by  

$$\P(\Bx_1+\Bf_1 \in H) \le \sup_v\P_{x_1,\dots,x_n}(|x_1v_1+\dots+x_nv_n-v|\le n^{-A}).$$ 

We again apply Theorem \ref{theorem:step1} to obtain a structural vector $\Bu$, and then use Theorem \ref{theorem:step2} to conclude that the probability for the existence of such $\Bu$ is negligible, completing the proof of Theorem \ref{theorem:discrete}.

The rest of the paper is organized as follows. In Section \ref{section:ILO} we introduce our key lemmas. Theorems \ref{theorem:step1} and \ref{theorem:step2} will be proven in Sections \ref{section:step1} and \ref{section:step2} respectively.

\section{The main tools for proving Theorem \ref{theorem:step1}}\label{section:ILO}
We need to study the concentration of $\sum_i x_iv_i$ in a small ball, where $\Bx=(x_1,\dots,x_n)$ is sampled uniformly from the set $\mathcal{S}$ of all $(-1,1)$ vectors of row-sum either $s-1$ or $s+1$. As customary, we first study a similar problem for $\Bx'$, a random vector whose components are i.i.d. copy of the Bernoulli variable $x$ defined in \eqref{eqn:Bernoulli}.

Let $V=\{v_1,\dots,v_n\}$ be a multiset in $\R^d$, where $d$ is a fixed integer. For $\beta>0$, we define the {\it small ball probability} as

$$\rho_{\beta}(V):=\sup_{v\in \R} \P_{\Bx'}\big(v_1x_1'+\dots+v_nx_n' \in B(v,\beta)\big),$$

where by $B(v,\beta)$ we denote the closed disk of radius $\beta$ centered at $v$ in $\R^d$. 

A well-known result of Erd\H{o}s \cite{E} and Littlewood-Offord \cite{LO}  asserts that if  $v_i$ are real numbers of magnitude $|v_i|\ge \beta$, then 

$$\rho_{\beta}(V)= O(n^{-1/2}).$$ 

This remarkable inequality has generated an impressive way of research. We refer the reader to \cite{H,Kle,NgV-optimal,TVinverse} and the references therein for further discussion regarding these developments.

In the reverse direction, we would like to find the underlying reason as to why the small ball probability is large (say, polynomial in $n$). 

Typical examples of $V$, where $\rho_{\beta}$ is large, involve generalized arithmetic progressions introduced in the previous section.
 
\begin{example}\label{example:linear}
Let $Q= \{\sum_{i=1}^r k_ig_i | -K_i \le k_i \le K_i\}$ be a proper symmetric GAP of rank $r=O(1)$ and size $N=n^{O(1)}$ in $\R^d$. Assume that for each $v_i$ there exists $q_i\in Q$ such that $\|v_i-q\|\le \delta$. Then, because the random sum $\sum_i q_ix_i'$ takes value in the GAP $nQ:=\{\sum_{i=1}^r k_ig_i| -nK_i \le k_i \le nK_i\}$, and because $|nQ| \le n^r N=n^{O(1)}$, the pigeon-hole principle implies that $\sum_i q_ix_i$ takes some value in $nQ$ with probability $n^{-O(1)}$. Thus we have

\begin{equation}\label{bound2} 
\rho_{n\delta}(V)   = n^{-O(1)}.
\end{equation}
\end{example}

The above example shows that if $v_i$ are {\it close} to a $GAP$ of rank $O(1)$ and size $n^{O(1)}$ in $\R^d$, then $V$ has large small ball probability. It was shown by Tao and the second author in \cite{TVbull,TVcir,TVinverse,TVcomp}, and by the current authors in \cite{NgV-optimal} that these are essentially the only examples of large small ball probability. We present here a somewhat optimal version.

We say that a vector $v$ is {\it $\delta$-close} to a vector $q$ if $\|v-q\|\le \delta$. We say that $v$ is $\delta$-close to a set $Q$ if there exists $q\in Q$ such that $v$ is $\delta$-close to $q$.

\begin{theorem}[Continuous Inverse Littlewood-Offord theorem for Bernoulli distribution]\cite[Theorem 2.9]{NgV-optimal}\label{theorem:ILO}
Let $0 <\alpha < 1/2; 0 < C$ be constants. Let $ \beta >0$ be a parameter that may depend on $n$. Suppose that $V=\{v_1,\dots,v_n\}$ is a multi-subset of $\R^d$ such that $\sum_{i=1}^n\|v_i\|^2=1$ and that $V$ has large small ball probability

$$\rho:= \rho_{\beta}(V)\ge n^{-C},$$

where in the definition of $\rho_{\beta}$ we assume $x_1',\dots,x_n'$ to be i.i.d. copies of the Bernoulli random variable $x$ defined in \eqref{eqn:Bernoulli}. Then for any number $n^\alpha \le n'  \le n$, there exists a proper symmetric GAP $Q=\{\sum_{i=1}^r k_ig_i : |k_i|\le K_i \}$ such that the following holds.

\begin{itemize}

\item (Full dimension) There exists $\sqrt{\frac{n'}{\log n}} \ll k \ll \sqrt{n'}$ such that the dilate $P:=(\beta/k)^{-1}\cdot Q$ contains the discrete hypercube $\{0,1\}^d $. Furthermore $P$ is an integral set, $P\subset \Z^d$.

\vskip .1in

\item (Approximation) At least $n-n'$ elements of $V$ (counting multiplicity) are $O(\frac{\beta}{k})$-close to $Q$.

\vskip .1in

\item (Small rank and cardinality) $Q$ has constant rank $d \le r=O(1)$, and small cardinality

$$|Q| =\max\Big(1,O_{\alpha,d,C}(\rho^{-1} n'^{(-r+d)/2})\Big).$$

\vskip .1in

\item (Small generators) There is a non-zero integer $p=O(\sqrt{n'})$ such that  all steps $g_i$ of $Q$ have the form  $g_i=(g_{i1},\dots,g_{id})$, where $g_{ij}=\beta \cdot \frac{p_{ij}} {p} $ with  $p_{ij} \in \Z$ and $p_{ij}=O(\beta^{-1} \sqrt {n'}).$

\end{itemize}
\end{theorem}

We note that \cite[Theorem 2.9]{NgV-optimal} was originally stated for more general distribution of the $x_i'$. Another slight difference is that we require $P$ to be a subset of $\Z^d$ here. However, this additional fact is not new as it has been explicitly verified in the proof of Theorem 2.9 (see the last part of \cite[Section 6]{NgV-optimal}).

\begin{remark}\label{remark:ILO}
As noticed in \cite[Corollary 2.10]{NgV-optimal}, the above theorem implies that if we use a coarser structure (which $O(\beta)$-approximates the $v_i$ rather than $O(\beta/k)$-approximates as stated in Theorem \ref{theorem:ILO}), then we can obtain a bound of at most $\max(O(\rho^{-1}/\sqrt{n'}),1)$ in the size of $Q$. As it turned out, the saving factor $1/\sqrt{n'}$ here plays a crucial role in any applications of Theorem \ref{theorem:ILO} in the literature. 
\end{remark}

From now on we will be mainly working with $\R^2$ (equivalently, $\C$). Our method naturally extends to $\R^d$ for any fixed $d$ but we do not attempt to do so here. To prove Theorem \ref{theorem:step1} we need to modify our notion of concentration probability as follows. Let $V=\{v_1,\dots,v_n\}$ be a multiset in $\R^2$. For any $\beta>0$, we define 

$$\rho_{\beta}^\ast(V):=\sup_{v\in \R^2} \P_{\Bx}\big(v_1x_1+\dots+v_nx_n \in B(v,\beta)\big),$$

where the probability is taken uniformly over all $(-1,1)$ vectors $\Bx=(x_1,\dots, x_n)$ of given entry sum $\bar{s}$, where $|\bar{s}|\le (1-\ep)n$. (In later application we will set $\bar{s}$ to be either $s-1$ or $s+1$.) 

By definition, $\rho^\ast$ is invariant under translation. One observes that for any $\beta$ and $V$ we have

\begin{equation}\label{eqn:rho:relation}
\rho_\beta(V)=\Omega(\rho_\beta^\ast(V)/\sqrt{n}).
\end{equation}

This relation suggests that if $\rho^\ast:=\rho_\beta^\ast(V)$ is large, then Theorem \ref{theorem:ILO} (more precisely, Remark \ref{remark:ILO}) implies that all the $v_i$ can be approximated by a $GAP$ $Q$ of size $O((\rho^\ast)^{-1}\sqrt{n}/\sqrt{n'})$. This bound, unfortunately, falls short for any application as the saving factor $\sqrt{n}/\sqrt{n'}$ here is  greater than 1 (we refer the reader to Remark \ref{remark:non-degenerate} of Section \ref{section:step2} for more explanation). 

The above discussion shows that  a sole application of \eqref{eqn:rho:relation} is not enough to obtain a useful inverse result regarding $\rho^\ast$. In the following result, by using the extra translation invariance property of $\rho^\ast$, we provide a more economical inverse result.

\begin{theorem}[Inverse Littlewood-Offord result with respect to $\rho^\ast$]\label{theorem:ILO:ast}
Suppose that $V=\{v_1,\dots,v_n\}$ is a multi-subset of $\R^2$ such that $\sum_{i=1}^n\|v_i\|^2=1$ and that

$$\rho^\ast:= \rho_{\beta}^\ast (V)\ge n^{-C}$$

for some $\beta=O(n^{-21C-12})$. Then for any number $n^{\alpha} \le n'  \le n$ there exists a proper GAP $Q^\ast=\{g_0+\sum_{i=1}^r k_ig_i : |k_i|\le K_i \}$ such that

\begin{itemize}

\item At least $n-n'$ elements of $V$ are $\beta n^{5C+3}$-close to $Q^\ast$.

\vskip .1in

\item $Q^\ast$ has small rank $r=O (1)$, and small cardinality

$$|Q| =\max\Big(1,O_{\alpha,C}((\rho^\ast)^{-1}\sqrt{n}/n')\Big).$$

\vskip .1in

\item There is a non-zero integer $p=O(\sqrt{n'})$ such that  all steps $g_i=(g_{i1},g_{i2}), 0\le i\le r$ of $Q^\ast$ have the form  $g_{ij}=\beta \cdot \frac{p_{ij}} {p} $ with  $p_{ij} \in \Z$ and $p_{ij}=O(\beta^{-1} \sqrt {n'}).$

\end{itemize}

\end{theorem}

Note that the approximation in this case is not as fine as in Theorem \ref{theorem:ILO}(or as in Remark \ref{remark:ILO}) and the structure $Q^\ast$ is not necessarily symmetric. On the other hand, the size of $Q^\ast$ is bounded by $O((\rho^\ast)^{-1}\sqrt{n}/n')$, which is considerably smaller than $O((\rho^\ast)^{-1}\sqrt{n}/\sqrt{n'})$ obtained by \eqref{eqn:rho:relation}.

Before proving Theorem \ref{theorem:ILO:ast}, let us provide a useful fact whose proof is simple and hence omitted. 

\begin{fact}\label{fact:1} Assume that $P=\{k_1g_1+\dots+k_rg_r | -K_i \le k_i \le K_i\}$ is a proper symmetric GAP which contains $w_1,\dots, w_r$, where each $w_i$ can be written as $k_{i1}g_1+\dots+k_{ir}g_r, k_{ij}\in \Z, |k_{ij}|\le K_i$.

\begin{enumerate}[(i)] 

\item Assume that the vectors $\Bk_i=(k_{i1},\dots,k_{ir}), 1\le i\le r,$ have full rank in $\R^{r}$. Then we can express each generator $g_i$ as $g_i=y_{i1} w_1+\dots + y_{ir} w_r$, where $y_{ij}$ are rational numbers of the form $p/q$ with $|p|,|q| =O_r(|P|^r)$.

\vskip .2in

\item Assume that $\Bk_r$ belongs to the space spanned by $\Bk_1,\dots,\Bk_{r-1}$, then we can write $\Bk_r$ as $\Bk_r=y_1\Bk_1+\dots+y_{r-1}\Bk_{r-1}$, where $y_i$ are rational numbers of the form $p/q$ with $|p|,|q|=O_r(|P|^r)$.

\end{enumerate}
\end{fact}

We now proceed to justify the main result of this section.

\begin{proof}(of Theorem \ref{theorem:ILO:ast}) Define a new set $U\subset \R^3$ as

$$U=\{u_1,\dots,u_n\}:=\Big\{\frac{1}{2}\cdot(v_1,\frac{1}{\sqrt{n}}),\dots,\frac{1}{2}\cdot (v_n,\frac{1}{\sqrt{n}})\Big\}.$$ 

By definition, we have $\sum_i \|u_i\|^2=1$ and $\rho_\beta^\ast(V)=\rho_{\beta/2}^\ast(U)$. Thus, by \eqref{eqn:rho:relation}

$$\rho_{\beta/2}(U)=\Omega(\rho_{\beta/2}^\ast(U)/\sqrt{n})=\Omega(\rho^\ast_\beta(V)/\sqrt{n})=\Omega(n^{-C-1/2}).$$

We apply Theorem \ref{theorem:ILO} to $U$ to obtain two GAPs $Q$ and $P=(\beta/2)^{-1}k\cdot Q$ respectively. First, observe that  if the rank $r$ of $Q$ (and $P$) is at least 5, then

$$|Q|=O((\rho^\ast)^{-1}\sqrt{n}/n'^{(r-3)/2})=O((\rho^\ast)^{-1}\sqrt{n}/n'),$$ 

and so we are done by letting $Q^\ast$ be the GAP generated by the first two coordinates of the generators of $Q$. Note that $g_0=0$ because $Q$ is homogeneous. Also, we obtained a very good approximation (of order $O(\beta/k)$) in this case.

Next we observe that $r$ cannot be 3.  Assume otherwise that $P=\{\sum_{i=1}^3 k_ig_i:|k_i|\le K_i\}$, where $g_i=(g_{i1},g_{i2},g_{i3})\in \Z^3$ are the generators of $P$. Because $P\subset \Z^3$ and it contains $(1,0,0),(0,1,0)$ and $(0,0,1)$, by Fact \ref{fact:1} (i) the generators $g_i$ must have the form $(g_{i1},g_{i2},g_{i3})$ where $|g_{ij}|$ are bounded by $O(|P|^3)$. But $P$ has size $O(\rho_{\beta/2}^{-1}(U))=O(n^{C+1/2})$, thus $|g_{ij}|=O(n^{3C+3/2})$. As a consequence, all of the elements of $P$ must have norm at most $O(n^{4C+2})$. However, this is impossible because as one of the elements of $P$ is $O(1)$-close to an element of $(\beta/2k)^{-1}\cdot U$, its second coordinate must be of order at least  $\frac{\beta^{-1}k}{\sqrt{n}}$, which is greater than $n^{4C+2}$ by the assumption of $\beta$ of being sufficiently small. 

We now consider the case $r=4$, $P=\{\sum_{i=1}^4 k_ig_i:|k_i|\le K_i\}$, where $g_i=(g_{i1},g_{i2},g_{i3})\in \Z^3$. 
Let $(w_1,l),\dots,(w_{n-n'},l)$ be the elements of $P$ which are $O(1)$-close to $n-n'$ elements of the dilated set $(\beta/2k)^{-1} \cdot U$. Apparently $l=\Theta(\beta^{-1}k/\sqrt{n})$. We next consider two cases.

{\bf Case 1.} If all $\|w_i\|$ are smaller that $n^{4C+2}$, then we would be done because in this case the order of all $\|u_i\|$ is at most $O((\beta/2k) n^{4C+2})$, which is bounded by $\beta n^{4C+2}$.
 
{\bf Case 2.} Assume otherwise that, say $\|w_1\| \ge n^{4C+2}$. Consider the following elements of $P$, $\Bb_1:=(1,0,0), \Bb_2:=(0,1,0), \Bb_3:=(0,0,1)$ and $\Bb_4:=(w_1,l)$. Because $\|w_1\|$ is greater than $n^{4C+2}$, one checks that the condition of Fact \ref{fact:1} (ii) does not hold for $\Bb_1,\Bb_2,\Bb_3$ and $\Bb_4$. We thus apply Fact \ref{fact:1} (i) to conclude that each $g_i$ can be expressed as in the form $c_{i1}\Bb_1+c_{i2}\Bb_2+c_{i3}\Bb_3+c_{i4}\Bb_4$, where $c_{ij}=p/q$ and $|p|,|q|=O(n^{4C+2})$. 

Next, consider any $\Bb=(w_{i_0},l)$ from the set  $\{(w_1,l),\dots,(w_{n-n'},l)\}$. There exist $k_1,k_2,k_3,k_4\in \Z, |k_i|\le K_i$, such that $\Bb=k_1g_1+k_2g_2+k_3g_3+k_4g_4$, and so

\begin{align*}
\Bb&=(k_1c_{11}+k_2c_{21}+k_3c_{31}+k_4c_{41})\Bb_1+ (k_1c_{12}+x_2k_{22}+x_3c_{32}+k_4c_{42})\Bb_2 \\
&+ (k_1c_{13}+k_2c_{23}+k_3c_{33}+k_4c_{43})\Bb_3+ (k_1c_{14}+k_2c_{24}+k_3c_{34}+k_4c_{44})\Bb_4.
\end{align*}

Notice that $l=\Theta(\beta^{-1}k/\sqrt{n}) \ge n^{21C+11}$, meanwhile $|k_1c_{13}+k_2c_{23}+k_3c_{33}+k_4c_{43}| =O(n^{5C+5/2})$ and $|k_1c_{14}+k_2c_{24}+k_3c_{34}+k_4c_{44}|= \Theta( n^{-16C-8})$ as $c_{ij}$ are rational numbers whose denominators are bounded by $O(n^{4C+4})$ and $k_1c_{14}+k_2c_{24}+k_3c_{34}+k_4c_{44}$ cannot be zero. We conclude that the coefficients of $\Bb_3$ and $\Bb_4$ must be 0 and 1 respectively,

It thus follows that, by considering the first two coordinates of $\Bb_1$ and $\Bb_2$, 

\begin{align*}
\|w_{i_0}-w_1\|^2&= \big((k_1c_{11}+k_2c_{21}+k_3c_{31}+k_4c_{41})^2+  (k_1c_{12}+k_2c_{22}+k_3c_{32}+k_4c_{42})^2\big)^{1/2}\\
& = O(n^{5C+5/2})<n^{5C+3}.
\end{align*}

Combining Case 1 and 2, we infer that if $r=4$ then all but $n'$ elements of $V$ are $\beta n^{5C+3}$-close to a common point. To complete the proof, we just simply set $g_0=\beta n^{5C+3} \cdot p$ be this approximated point where $p$ is a complex number of integral coordinates and $|p|\le \beta^{-1} n^{-5C-3}$. We set other generators  to be zero.

\end{proof}

We now deduce an important corollary of Theorem \ref{theorem:ILO:ast} which, similarly to the result of Erd\H{o}s and Littlewood-Offord, states that as long as the multi-set $V$ is not too degenerated (for a given $\beta$), its concentration probability $\rho^\ast$ must be small.

\begin{corollary}\label{cor:LO}
Let $0<\alpha<1/2$ be a positive constant and let $n'$ be a number satisfying $n^{1/2+\alpha}<n'<n$. Assume that $\beta \le n^{-24}$ and $V$ is a multi-set in $\R^2$ so that any of its $n-n'$ elements cannot be $\beta n^{6}$-close to a common point. Then we have 

$$\rho_\beta^\ast(V) = O(\sqrt{n}/n').$$
\end{corollary}

\begin{proof}(of Corollary \ref{cor:LO})
Assume otherwise that $\rho_\beta^\ast(V) \ge C \sqrt{n}/n'$ for some large constant $C$ to be chosen. So 

$$\rho^\ast(V)\ge Cn^{-1/2}.$$ 

We next apply Theorem \ref{theorem:ILO:ast} to $V$ to obtain a GAP $Q^\ast$ which is $\beta n^{11/2}$ to all but $n-n'$ elements of $V$. Notice that because there are no more than $n-n'-1$ elements of $V$ that are $\beta n^{6}$-close to one common point, $Q^\ast$ must have size at least 2. On the other hand, from the conclusion of Theorem \ref{theorem:ILO:ast}, assuming that $C$ is sufficiently large depending on $\alpha$, the size of $Q^\ast$ is bounded by

$$|Q^\ast|=\max(1,O_\alpha((\rho^\ast)^{-1}\sqrt{n}/n')=\max(1,O_\alpha(\frac{1}{C}))=1.$$

This contradiction completes the proof of our corollary.

\end{proof}

\section{Proof of Theorem \ref{theorem:step1}}\label{section:step1}

We will invoke Theorem \ref{theorem:ILO:ast}. Define a radius sequence $(\beta_k)_0^\infty$ where $\beta_0:=n^{-A}$ and 

$$\beta_{i+1}=n^{5B+5+\gamma} \beta_i.$$

Let $V$ be the multi-set of $v_1,\dots,v_n$. Then the assumption of Theorem \ref{theorem:step1} becomes

$${\rho_{\beta_0}}^\ast(V)\ge n^{-B}.$$ 

with either $\bar{s}=s-1$ or $\bar{s}=s+1$.

Next, because the increasing sequence $\rho_{\beta_i}^\ast(V)$ is bounded from above by 1, by pigeonhole principle there exists $0\le k_0\le 2B/\alpha$  such that 

$$\rho_{\beta_{k_0+1}}^\ast(V)\le n^{\alpha/2}\rho_{\beta_{k_0}}^\ast(V).$$

As $A$ was chosen to be sufficiently large, one has $\beta_{k_0}\le n^{-A/2}$. We next apply Theorem \ref{theorem:ILO:ast} to $V$ with $n'=n^{1/2+\alpha}$ and $\beta=\beta_{k_0}$ to obtain a GAP $Q^\ast=\{g_0+\sum_{i=1}^r k_ig_i,|k_i|\le K_i\}$ for which the following holds.

\begin{itemize}

\item $Q^\ast$ has small rank $r=O (1)$, and small cardinality

$$|Q^\ast| =\max\Big(1,O_{\alpha,B}\big((\rho_{\beta_{k_0}}^\ast(V))^{-1}/n^\alpha\big)\Big).$$

\vskip .1in

\item There are $n_0:=n-n^{1/2+\alpha}$ elements $v_{i_1},\dots, v_{i_{n_0}}$ of $V$ which are $O(\beta_{k_0} n^{5B+3})$-close to $n-n^{1/2+\alpha}$ elements $u_1,\dots, u_{n_0}$ of $Q^\ast$.

\vskip .1in

\item There is a non-zero integer $p=O(\sqrt{n^{1/2+\alpha}})$ such that  all steps $g_i=(g_{i1},g_{i2}), 0\le i\le r$ of $Q^\ast$ have the form  $g_{ij}=\beta_{k_0}^{-1} p_{ij}/p $ with  $p_{ij} \in \Z$ and $p_{ij}=O(\beta_{k_0}^{-1} \sqrt {n^{1/2+\alpha}})$. In particular, all the components of the elements of $Q^\ast$ have the form $p/q$  where $|p|,|q| \le n^{A+1}$.

\end{itemize}

Next, for each $v$ of the remaining $n^{1/2+\alpha}$ exceptional elements of $V$ (which are not close to any element of $Q^\ast$), we trivially approximate it by a complex number $v$ whose components are rational numbers of the form $p/q$ with $|q|\le n^{A+1}$ such that $|u-v|\le \beta_{k_0} n^{5B+3}$.

By the approximation we infer that

$$\|\Bu-\Bv\|=(\sum_i |u_i-v_i|^2)^{1/2} \le \beta_{k_0} n^{5B+7/2}.$$ 

Taking into account that $|f_{ij}|\le n^\gamma$, we thus have

\begin{align*}
 \rho_{\beta_{k_0}}^\ast(V) &\le \rho_{\beta_{k_0}+ \beta_{k_0}n^{5B+7/2+\gamma}}^\ast(U)\le \rho_{\beta_{k_0}n^{5B+4+\gamma}}^\ast(U)\\
&\le \rho_{\beta_{k_0}+ \beta_{k_0}n^{5B+4+\gamma}}^\ast(V) \le \rho_{\beta_{k_0}n^{5B+5+\gamma}}^\ast(V) = \rho_{\beta_{k_0+1}}^\ast (V),
\end{align*}

where $U$ is the multi-set $\{u_1,\dots,u_n\}$. 

From the estimate above, as $\rho_{\beta_{k_0+1}}^\ast (V) \le n^{\alpha/2} \rho_{\beta_{k_0}}^\ast (V)$, it is implied that 

$$\rho_{\beta_{k_0}n^{5B+4+\gamma}}^\ast(U) \le n^{\alpha/2} \rho_{\beta_{k_0}}^\ast(V).$$ 

So the size of $Q^\ast$ is bounded by 

$$|Q^\ast|=\max\Big(1,O\big((\rho_{\beta_{k_0}n^{5B+4+\gamma}}^\ast(U))^{-1}/n^{\alpha/2}\big)\Big).$$

In summary, we have obtained a vector $\Bu=(u_1,\dots,u_{n})$ which satisfies the following properties.

\begin{itemize}
\item $\|\Bu\|\asymp 1$, and because $\langle \Bv, \row_i(X)\rangle =0$ for any row $\row_i$ of $X$ of index $i\ge 2$, we also have $|\langle \Bu,\row_i(X)\rangle| \le \beta_{k_0}n^{5B+4+\gamma} $.
\vskip .2in
\item There exists a generalized arithmetic progression $Q^\ast$ of rank $O_{B,\alpha}(1)$ and size  $|Q^\ast|=\max\Big(1,O\big((\rho_{\beta_{k_0} n^{5B+4+\gamma}}^\ast(U))^{-1}/n^{\alpha/2}\big)\Big)$ that contains at least $n-n^{1/2+\alpha}$ complex numbers $u_i$.
\vskip .2in
\item All the components of $u_i$ and of the generators of $Q^\ast$ are rational numbers of the form $p/q$, where $|p|,|q| \le n^{A+1}$.
\end{itemize}

This completes the proof of Theorem \ref{theorem:step1}.

\section{Proof of Theorem \ref{theorem:step2}}\label{section:step2}

By applying Theorem \ref{theorem:step1}, we obtain a structural vector $\Bu$ which satisfies all the described properties. Because the number of $\beta$ is bounded by a constant, it is enough to verify Theorem \ref{theorem:step2} for one such $\beta$. By paying a factor of $n$ in probability, we assume that $|\langle \Bu,\row_i(X)\rangle| \le \beta n^{5B+4+\gamma} $ for the last $n-1$ rows of $X$.

Set $\beta':=\beta n^{5B+4+\gamma}$. We will consider two cases depending on the structure of $\Bu$. 

\subsection{Degenerate $\Bu$} We first consider the probability $\P_{\bf major}$ of the event $|\langle \row_i,\Bu\rangle |\le \beta', 2\le i\le n $, for which there are $n_0:=n-n^{1/2+\alpha}$ complex numbers $u_i$ which can be $\beta' n^4$-approximated by a common point $u_0' \in \beta' n^4 \cdot \Z^2$. 

By paying a factor $\binom{n}{n_0}$ in probability, we may assume that this point approximates the first $n_0$ complex numbers $u_1,\dots, u_{n_0}$. Thus, by approximating the remaining $u_i$ by $u_i'\in \beta'n^4 \cdot \Z^2$  such that $|u_i-u_i'|\le \beta'n^4$, the events $|\langle \row_i,\Bu \rangle |\le \beta'$  belongs to the event $|\langle \row_i,\Bu' \rangle |\le \beta' n^5$, where $\Bu'=(u_1',\dots,u_1',u_{n_0+1}',\dots,u_n')$ and $\|\Bu'\|\asymp 1 $.

Let $X_{(n-1)\times n}$ be the matrix generated by the last $n-1$ rows of $X$, and let $X'$ be the $n-1$ by $n-n_0$ matrix obtained from $X_{(n-1)\times n}$ by joining its first $n_0$ columns, 

$$X'= \left[\col_1(X_{(n-1)\times n})+\dots+\col_{n_0}(X_{(n-1)\times n}),\col_{n_0+1}(X_{(n-1)\times n}), \dots ,\col_n(X_{(n-1)\times n})\right].$$

By definition, the row vectors of $X'$ satisfy $|\langle \row_i(X'),\Bu'_{\tr} \rangle|\le \beta' n^5$ where  $\Bu'_{\tr}:=(u_1',u_{n_0+1}',\dots,u_n')$. It also follows from definition that the $i$-th row of $X'$ has the form $\row_i(X')=\Bx' + \Bf'$, where $\Bf'=(f_{i1}+\dots+f_{in_0},f_{i(n_{0}+1)},\dots,f_{in})$ and $\Bx'=(x_1+\dots+x_{n_0},x_{n_0+1},\dots,x_n):=(x_1',\dots,x_{n-n_0}')$. 

As $\Bx$ is sampled uniformly from $\mathcal{S}$, the set of all $(-1,1)$ vectors of entry-sum either $s-1$ or $s+1$, $\Bx'$ is a random vector chosen from {\it type 1} or {\it type 2} defined below.

{\bf Type 1.} (row-sum $s+1$)

$$ 
\P(x_1'=k)=\frac{\binom{n_0}{(n_0+k)/2}\binom{n-n_0}{(n-n_0+s+1-k)/2}}{\binom{n}{n/2+(s-1)/2}+\binom{n}{n/2+(s+1)/2}} 
$$

for all $k$ such that $k+n_0$ is even; and $(x_2',\dots,x_{n-n_0}')$ are chosen uniformly from all $(-1,1)$ vectors of row-sum $s+1-x_1'$.

{\bf Type 2.} (row-sum $s-1$)

$$ 
\P(x_1'=k)=\frac{\binom{n_0}{(n_0+k)/2}\binom{n-n_0}{(n-n_0+s-1-k)/2}}{\binom{n}{n/2+(s-1)/2}+\binom{n}{n/2+(s+1)/2}}
$$

for all $k$ such that $k+n_0$ is even; and $(x_2',\dots,x_{n-n_0}')$ are chosen uniformly from all $(-1,1)$ vectors of row-sum $s-1-x_1'$.

It is clear that  

$$\P(\Bx' \in \mbox{type 1})= \frac{\binom{n}{n/2+(s+1)/2}}{\binom{n}{n/2+(s-1)/2}+\binom{n}{n/2+(s+1)/2}}$$

and 

$$\P(\Bx' \in \mbox{type 2})= \frac{\binom{n}{n/2+(s-1)/2}}{\binom{n}{n/2+(s-1)/2}+\binom{n}{n/2+(s+1)/2}}.$$

Observe that as $|s|\le (1-\ep)n$, these two probabilities are comparable, each of which can be bounded crudely from below by $(1-\ep)/4$.

We next apply the following result.

\begin{claim}\label{claim:1} Let $\ep<1/4$ be a fixed constant. Let $\Bu'_{\tr}=(u_1',u_{n_0+1}',\dots,u_n')$ be a vector in which the components of each complex $u_i'$ is of the form $\beta' n^4 \cdot \Z$ and such that  $n_0|u_1'|^2+|u_{n_0+1}'|^2+\dots+ |u_n'|^2 \asymp 1$. Then, as $n$ is sufficiently large and $\Bf$ is a fixed vector, one has 

$$\P_{\Bx'}(|\langle \Bx' + \Bf',\Bu'_\tr \rangle|\le \beta'n^5)\le 1-(1-\ep)/8.$$   
\end{claim}

\begin{proof}(of Claim \ref{claim:1}) We will consider two main cases below. 

{\bf (i)} We first assume that there exists $1<i_0<j_0 $ such that $|u_{i_0}'- u_{j_0}'|\ge \beta'n^5$. Without loss of generality, assume that $i_0=n-1$ and $j_0=n$. It follows from the distribution of $\Bx'$ that the event of having exactly one $-1$ among the last two components of $\Bx'$ happens with probability at least $(1-\ep)/4$ asymptotically. Within this event, observe that for any tuple $(x_1',\dots,x_{n-n_0-2}')$, either $\Bx=(x_1',\dots,x_{n-n_0-2}',-1,1)$ or $\Bx=(x_1',\dots,x_{n-n_0-2}',1,-1)$ does not satisfy $|\langle \Bx',\Bu'_\tr \rangle + \langle \Bf',\Bu'_\tr \rangle|\le \beta'n^5$. Thus we have 

$$\P_{\Bx'}(|\langle \Bx',\Bu'_\tr \rangle + \langle \Bf',\Bu'_\tr \rangle|\le \beta'n^5)\le 1-(1-\ep)/8.$$

{\bf (ii)} Assume otherwise that there exists $u'$ such that all $|u'-u_{n_0+1}'|, \dots, |u'-u_{n}'|$ are bounded by $\beta' n^5$. In this case, the inequality $|\langle \Bx',\Bu'_\tr \rangle + \langle \Bf',\Bu'_\tr \rangle|\le \beta' n^5$ implies that 

\begin{equation}\label{eqn:uniform}
|x_1'(u_1'-u')+u'(x_1'+x_2'+\dots+x_{n-n_0}')+\langle \Bf',\Bu'_\tr \rangle | \le \beta' n^6.
\end{equation}

We next consider the subcase $|u_1'-u'|\ge \beta' n^8$. If $x_1'+x_2'+\dots+x_{n-n_0}'=s+1$, then \eqref{eqn:uniform} implies that $x_1'$ belongs to the interval $\big[(-u'(s+1)-\langle \Bf',\Bu'_\tr \rangle)(\beta' n^8)^{-1} -1/n^2,(-u'(s+1)-\langle \Bf',\Bu'_\tr \rangle)(\beta' n^8)^{-1}+1/n^2\big]$. However, because this interval has length $2/n^2$, and so this probability is clearly bounded by $\sup_k \P(x_1'=k)$, which is clearly smaller than $1-(1-\ep)/4$. We argue similarly for the case $x_1'+x_2'+\dots+x_{n-n_0}'=s-1$.

For the remaining  subcase $|u_1'-u'| \le \beta' n^8$, as $A$ was chosen to be large enough, we have $|u_1'-u'|\le n^{-2}$. Next, because 
$\|\Bu'_\tr\|^2=n_0|u_0'|^2 + (n-n_0)|u'|^2 \asymp 1$, we infer that $|u'|\asymp 1/\sqrt{n}$. It then follows that 

\begin{equation}\label{eqn:uniform'}
|u'(x_1'+\dots+x_{n-n_0}')+\langle \Bf',\Bu'_\tr \rangle| \le \beta' n^9.
\end{equation}

However, as $x_1'+\dots+x_{n-n_0}'$ takes value $s+1$ and $s-1$ each with probability at least $(1-\ep)/4$, the equation \eqref{eqn:uniform'} above holds with probability at most $1-(1-\ep)/4$.

\end{proof}

Now we estimate $\P_{\bf major}$. As the event $|\langle \row_i(X),\Bu' \rangle|\le \beta' n^5$ is controlled by $|\langle \row_i(X'),\Bu'_\tr \rangle|\le \beta' n^5$, and By Claim \ref{claim:1} the later holds with probability $(7+\ep)/8$, it follows that the probability that $|\langle \row_i(X),\Bu' \rangle|\le \beta' n^5$ for all $2\le i\le n$ is bounded by $((7+\ep)/8)^{n-1}$. 

Additionally, an elementary computation implies that the number of structural vectors $\Bu'\in (\beta'n^4 \cdot \Z^2)^{n-n_0+1}$ satisfying $\|\Bu'\|\asymp 1$ is bounded by 

$$((\beta'n^4)^{-1})^{n-n_0+1} = O((n^A)^{n^{1/2+\ep}+1}) =O(n^{O_A(n^{1/2+\ep})}).$$

Putting together, we obtain the following bound for $\P_{\bf major}$

$$\P_{\bf major} = O(n^{O_A(n^{1/2+\ep})}) \binom{n}{n_0}\binom{n-1}{n-n_0-1} (\frac{7+\ep}{8})^{n-1}= (\frac{7+\ep}{8})^{(1-o(1))n}.$$

\vskip .2in

\begin{remark}
In the treatment above the fact that $\Bx'$ takes either type 1 or type 2 with comparable probability is crucial. The assumption of just one type would not be enough to estimate $\P_{\bf major}$ unless we had an additional assumption on $\Bu'$, say $u_1'+\dots+u_n'$ is nearly zero.     
\end{remark}

\subsection{Non-degenerate $\Bu$} We consider the probability $\P_{\bf minor}$ of the event that there exists a vector $\Bu$ for which $|\langle \row_i(X),\Bu \rangle|\le \beta', 2\le i$ and the following holds

\begin{itemize}
\item  $\|\Bu\|\asymp 1$ and there does not exist any $u$ which is $\beta'n^4$-close to all but $n^{1/2+\alpha}$ complex numbers $u_i$. Thus it follows from Corollary \ref{cor:LO} that 

$$\rho_{\beta'}^\ast(U)=O(n^{-\alpha}).$$
\vskip .2in
\item There exists a generalized arithmetic progression $Q^\ast$ of rank $O_{B,\alpha}(1)$ and size  $|Q^\ast|=\max\big(1,O(\rho_{\beta'}^\ast(U)^{-1}/n^{\alpha/2})\big)=O(\rho_{\beta'}^\ast(U)^{-1}/n^{\alpha/2})$ that contains at least $n-n^{1/2+\alpha}$ complex numbers $u_i$. (Here we used the estimate $\rho_{\beta'}^\ast(U)^{-1}= \Omega(n^\alpha)$ to eliminate the trivial constant 1 in the size estimate of $Q^\ast$.) 
\vskip .2in
\item All the components of $u_i$ and of the generators of the generalized arithmetic progression are rational numbers of the form $p/q$, where $|p|,|q| \le n^{A+1}$.
\end{itemize}

Let $0<\delta$ to be chosen (any $\delta <\alpha/3$ will suffice) . We divide the interval $[n^{-B},O_\alpha(n^{-\alpha/2})]$ into sub-intervals $[n^{-(k+1)\delta},n^{-k\delta}]$, where $\alpha/2\delta \le k\le B/\delta$.  For each $k$, let $\BG_k$ be the collection of $\Bu$'s such that $\rho_{\beta'}^\ast(U)\in [n^{-(k+1)\delta},n^{-k\delta}]$, and let $\P_k$ be the probability that $|\langle \row_i(X'),\Bu\rangle |\le \beta'$ for all $i$ and for one of $\Bu$ from $\BG_k$. 

We now bound the size of $\BG_k$. To do this, we first count the number of GAPs which may contain most of the $u_i$ of vectors $\Bu$ from $\BG_k$, and then count the number of $\Bu$'s whose $u_i$ are chosen from the determined structure. Recall that all components of the GAP generators are of the form $p/q$, where $|p|,|q|\le n^{A+1}$. Because each GAP has rank $O_{B,\alpha}(1)$ and size $O((\rho^\ast)^{-1}/n^{\alpha/2})=O(n^{\delta(k+1)}/n^{\alpha/2})$, the number of such GAPs is bounded by 

$$(n^{4A+4})^{O_{B,\alpha}(1)} (n^{\delta(k+1)}/n^{\alpha/2})^{O_{B,\alpha}(1)}= O(n^{O_{B,\alpha,\delta}(1)}).$$

After choosing a $Q^\ast$ of size $O(n^{\delta(k+1)}/n^{\alpha/2})$, the number of ways to choose $n-n^{1/2+\alpha}$ complex numbers $u_i$ as $Q^\ast$'s elements is

$$\binom{n}{n^{1/2+\alpha}}\binom{O(n^{\delta(k+1)}/n^{\alpha/2})}{n-n^{1/2+\alpha}}=O\big(n^{n^{1/2+\alpha}}(n^{\delta(k+1)}/n^{\alpha/2})^{n-n^{1/2+\alpha}}\big).$$

For the remaining $n^{1/2+\alpha}$ exceptional elements, there are $(n^{4A+4})^{n^{1/2+\alpha}}=O(n^{O_A(n^{1/2+\alpha})})$ ways to choose them. Putting these bounds together, we obtain the following bound for the number of $\Bu$ of $\BG_k$

$$|\BG_k|=O\big(n^{O_{A,B,\alpha,\delta}(n^{1/2+\alpha})}(n^{\delta(k+1)}/n^{\alpha/2})^{n-n^{1/2+\alpha}}\big).$$

Now, for a given $\Bu\in \BG_k$, the probability that  $|\langle \row_i(X),\Bu\rangle |\le \beta'$ for all $2\le i\le n$ is bounded by $(\rho_{\beta'}^\ast(\Bu))^{n-1}\le (n^{-\delta k})^{n-1}$. Thus we can estimate $\P_k$ as 

$$\P_k \le |\BG_k|(n^{-\delta k})^{n-1} =O\left(n^{O_{A,B,\alpha,\delta}(n^{1/2+\alpha})} (n^\delta)^{n}/(n^{\alpha/2})^{n-n^{1/2+\alpha}}\right)=o(n^{-\alpha n/6}),$$

provided that $\delta$ was chosen to be smaller than $\alpha/3$.

Summing over $k$, we thus obtain

$$\P_{\bf minor}=\sum_{k\le B/\delta} \P_k= o(n^{-\alpha n/6}).$$

\begin{remark}\label{remark:non-degenerate}
One observes that the saving factor $1/n^{\alpha/2}$ in the size of $Q^\ast$ plays a key role in our analysis here. This explains the necessity of Theorem \ref{theorem:ILO:ast}.
\end{remark}

\section{Concentration of distance}\label{section:Talagrand}

We now give a proof of Lemma \ref{lemma:Talagrand} basing on \cite{TVdet}. Let $P=(p_{ij})$ be the $n$ by $n$ orthogonal projection matrix from $\C^n$ to $V^\perp$. Thus $P$ is Hermitian and $P^2=P$. We first normalize $x_i'$ by setting $y_i':=x_i'-s/n$ and $f_i':=f_i+s/n$ for $1\le i\le n$. We then have $\E y_i'=0, \Var(y_i') =1-(s/n)^2$ and

\begin{align*}
d'^2 = \|P(\Bf+\Bx')\|^2&= \|P(\Bf'+\By')\|^2= \sum_{ij} p_{ij}(y_i'+f'_i)\overline{(y_j'+f'_j)}\\ 
&= \sum_{ij} p_{ij} y_i'y_j' + \sum_{ij}y_i'(p_{ij}\overline{f'_j}+p_{ji}f'_j) + \sum_{ij}p_{ij}f'_i\overline{f'_j}\\
&= \Tr(P) + \sum_{i\neq j} p_{ij}y_i'y_j' + \sum_{ij}y_i'(p_{ij}\overline{f'_j}+p_{ji}f'_j) +  \sum_{ij}p_{ji}y_i'f'_j + d_{\Bf'}^2\\
&:= (n-k) + d_{\Bf'}^2 + Y.
\end{align*}

It is clear that $\E Y=0$, thus

$$\E(d^2)=(n-k)+d_{\Bf'}^2.$$

Note that 

\begin{align*}
\E |Y|^2 &= \E | \sum_{i\neq j} p_{ij}y_i'y_j' + \sum_{ij}y_i'(p_{ij}\overline{f'_j}+p_{ji}f'_j)|^2 \\
&=\E |\sum_{i\neq j} p_{ij}y_i'y_j'|^2 +   \E |\sum_{ij}y_i'(p_{ij}\overline{f'_j}+p_{ji}f'_j)|^2 \\
&=(1-(s/n)^2)\Big[\sum_{i\neq j}|p_{ij}|^2 + \sum_i |\sum_j p_{ij}\overline{f'_j} +  \sum_{j}p_{ji}f'_j|^2 \Big]\\
&\le \sum_{i\neq j}|p_{ij}|^2 + 4\sum_i (\Re(\sum_{j}p_{ji}f'_j))^2\\
&\le \sum_{i\neq j}|p_{ij}|^2 + 4\sum_i |\sum_{j}p_{ji}f'_j|^2\\
&=\sum_{i\neq j}|p_{ij}|^2 + 4\sum_{j_1j_2} \sum_{i}p_{j_1i} \overline{p_{j_2i}} f'_{j_1}\overline{f'_{j_2}}\\
&=\sum_{i\neq j}|p_{ij}|^2 + 4\sum_{j_1j_2} p_{j_1j_2}f'_{j_1}\overline{f'_{j_2}} =\sum_{i\neq j}p_{ij}^2 + 4d_{\Bf'}^2.
\end{align*}

Next, because $\sum_i p_{ii}=(n-k)$, by Cauchy-Schwarz inequality 

$$\sum_i p_{ii}^2 \ge (n-k)^2/n.$$

Thus 

$$\sum_{i\neq j}|p_{ij}|^2 = \sum_{i,j}|p_{ij}|^2 -\sum_{i}p_{ii}^2 \le (n-k) - (n-k)^2/n \le \min(k,n-k).$$

It is implied that 

$$\E Y^2 \le \min(k,n-k) +4d_{\Bf'}^2.$$

Consider the event $d \ge \sqrt{n-k + d_{\Bf'}^2} + 3$. The probability of this event is bounded from above by 

\begin{align*}
&\P\big(d'^2 \ge n-k+d_{\Bf'}^2 + 6  \sqrt{n-k + d_{\Bf'}^2}\big)\\
&=P\big(Y\ge  6  \sqrt{n-k + d_{\Bf'}^2}\big) \\
&\le \P\big(Y^2 \ge 36  (n-k + d_{\Bf'}^2)\big)\\
&\le \frac{\E Y^2}{36(n-k+d_{\Bf'}^2)} \le \frac{1}{9}. 
\end{align*}

Similarly, consider the event  $d' \le \sqrt{n-k + d_{\Bf'}^2} - 3$. The probability of this event is bounded from above by 

\begin{align*}
&\P\big(d'^2 \le n-k+d_{\Bf'}^2 - 6  \sqrt{n-k + d_{\Bf'}^2}+9\big)\\
&=P\big(Y\le  -6  \sqrt{n-k + d_{\Bf'}^2}+9\big) \\
&\le \P\big(Y^2 \ge 36  (n-k + d_{\Bf'}^2) - 108  \sqrt{n-k + d_{\Bf'}^2}+81\big)\\
&\le \frac{\E Y^2}{36  (n-k + d_{\Bf'}^2) - 108  \sqrt{n-k + d_{\Bf'}^2}+81} \le \frac{1}{4} 
\end{align*}

provided that $k\le n-10$.

Thus the median $M$ of $d'$ satisfies $|M-\sqrt{n-k + d_{\Bf'}^2}|\le 3$.

Since the distance function is convex on $\{-1,1\}^n$ with Lipschitz constant 1. Talagrand's concentration inequality \cite{T} implies that for any $t$

$$\P(|d' -M|\ge t)\le 4\exp(-t^2/16).$$

Since $|M-\sqrt{n-k + d_{\Bf'}^2}|\le 3$, Lemma \ref{lemma:Talagrand} follows.

\end{document}